\documentclass[12pt]{article}
\usepackage[dvips]{graphicx}
\usepackage{amsmath}
\usepackage{amsthm}
\usepackage{amsfonts}
\usepackage{amssymb}

\newtheorem{thm}{Theorem}

\newtheorem{question}{Question}

\Alph{thm}

\theoremstyle{definition}

\newtheorem{rem}[thm]{Remark}

\newcommand{\Integer}{\mathbb{Z}}

\newcommand{\F}{{\mathrm F}}
\newcommand{\Aut}{\operatorname{Aut}}
\newcommand{\Inn}{\operatorname{Inn}}

\newcommand{\Isom}{\operatorname{Isom}}

\newcommand{\X}{{\langle x \rangle}}

\newcommand{\co}{\!:\!}
\newcommand{\CATZ}{{\rm CAT}(0)}

\newcommand{\subheading}[1]{\vskip 10pt \noindent {\bf #1}: }

\usepackage{ifthen}

\newcommand{\showcomments}{yes}

\newsavebox{\commentbox}
{\ifthenelse{\equal{\showcomments}{yes}}%
{\footnotemark
    \begin{lrbox}{\commentbox}
    \begin{minipage}[t]{1.25in}\raggedright\sffamily\upshape\tiny
    \footnotemark[\arabic{footnote}]}
{\begin{lrbox}{\commentbox}}}%
{\ifthenelse{\equal{\showcomments}{yes}}%
{\end{minipage}\end{lrbox}\marginpar{\usebox{\commentbox}}}
{\end{lrbox}}}

\begin{document}
\title{The automorphism group of the free group of rank two is a $\CATZ$ group}
\author{Adam Piggott\footnote{Department of Mathematics, Bucknell University, Lewisburg, PA 17837, USA }, Kim Ruane\footnote{Department of Mathematics, Tufts University, Medford, MA 02155, USA} and Genevieve S. Walsh\footnote{Department of Mathematics, Tufts University, Medford, MA 02155, USA.  This author was supported by N. S. F. grant DMS-0805908.}}
\date{\today}
\maketitle

\begin{abstract}
We prove that the automorphism group of the braid group on four
strands acts faithfully and geometrically on a $\CATZ$ 2-complex.
This implies that the automorphism group of the free group of rank
two acts faithfully and geometrically on a $\CATZ$ 2-complex, in
contrast to the situation for rank three and above.
\end{abstract}

\section{Introduction}

A \emph{CAT(0) metric space} is a proper complete geodesic metric
space in which each geodesic triangle with side lengths $a$, $b$ and
$c$ is ``at least as thin'' as the Euclidean triangle with side
lengths $a$, $b$ and $c$ (see \cite{MartinsBook} for details). We
say that a finitely generated group $G$ is a \emph{CAT(0) group} if
$G$ may be realized as a cocompact and properly discontinuous
subgroup of the isometry group of a $\CATZ$ metric space $X$.
Equivalently, $G$ is a $\CATZ$ group if there exists a $\CATZ$
metric space $X$ and a faithful geometric action of $G$ on $X$. It
is perhaps not standard to require that the group action be
faithful, a point which we address in Remark \ref{QuestionRemark}
below.

For each integer $n \geq 2$, we write $F_n$ for the free group of
rank $n$ and $B_n$ for the braid group on $n$ strands.

In \cite{AutomaticStructuresOnAutF2}, T. Brady exhibited a subgroup
$H \leq \Aut(F_2)$ of index 24 which acts faithfully and
geometrically on $\CATZ$ 2-complex.  In subsequent work
\cite{BradyArtinGroupsWithThreeGenerators}, the same author showed
that $B_4$ acts faithfully and geometrically on a $\CATZ$ 3-complex.
It follows that $\Inn(B_4)$ acts faithfully and geometrically on a
$\CATZ$ 2-complex $X_0$ (this fact is explained explicitly by Crisp
and Paoluzzi in
\cite{CrispPaoluzzi}). 
Now, $\Inn(B_n)$ 
has index two in $\Aut(B_n)$ \cite{AutomorphismsOfBraidGroups}, and
$\Aut(F_2)$ is isomorphic to $\Aut(B_4)$ \cite{B4IsLinear,
AutomorphismsOfBraidGroups}, thus the result in the title of this
paper is proved if we exhibit an extra isometry of $X_0$ which
extends the faithful geometric action of $\Inn(B_4)$ to a faithful
geometric action of $\Aut(B_4)$. We do this in
$\S$\ref{AutB4Section} below.

In the language of \cite{SimplicalNonpositiveCurvature}, $X_0$ is a systolic simplicial complex.  By \cite[Theorem 13.1]{SimplicalNonpositiveCurvature}, a group which acts simplicially, properly discontinuously and cocompactly on
such a space is biautomatic.  Since the action of $\Aut \F_2$ provided here is of this type, it follows that $\Aut \F_2$ is biautomatic.

Our results reinforce the striking contrast between those properties
enjoyed by $\Aut(F_2)$ and those enjoyed by the automorphism groups
of finitely generated free groups of higher ranks.  We can now say
that $\Aut(F_2)$ is a $\CATZ$ group, a biautomatic group and it has a faithful linear
representation \cite{BraidsAndAutFTwo, B4IsLinear}; while
$\Aut(F_n)$ is not a $\CATZ$ group \cite{GerstenAutFNotCAT0}, nor a
biautomatic group \cite{AutFnNotBiautomaticFornAtLeast3} and it does
not have a faithful linear representation
\cite{AutFNotLinearForNGEQ3} whenever $n \geq 3$.

We regard the $\CATZ$ 2-complex $X_0$ as a geometric companion to
Auter Space (of rank two) \cite{AuterSpace}, a topological
construction equipped with a group action by $\Aut(F_2)$.

Let $W_3$ denote the universal Coxeter group of rank $3$---that is,
$W_3$ is the free product of $3$ copies of the group of order two.
 Since $\Aut(F_2)$ is isomorphic to $\Aut(W_3)$ (see Remark
\ref{IsomorphismRemark} below), we also learn that $\Aut(W_3)$ is a
$\CATZ$ group.


\begin{rem} \label{QuestionRemark} As pointed out in the opening paragraph, our definition of a $\CATZ$ group
is perhaps not standard because of the requirement that the group
action be faithful.  We note that such a requirement is redundant
when giving an analogous definition of a word hyperbolic group. This
follows from the fact that word hyperbolicity is an invariant of the
quasi-isometry class of a group.
In contrast, the $\CATZ$ property is not an invariant of the
quasi-isometry class of a group.
Examples are known of two quasi-isometric groups, one of which is
$\CATZ$, and the other of which is not. Examples of this type may be
constructed using the fundamental groups of graph manifolds
\cite{KapovichLeeb} and the fundamental groups of Seifert fibre
spaces \cite[p.258]{MartinsBook}\cite{BridsonAlonso}. So the
adjective `faithful' is not so easily discarded in our definition of
a $\CATZ$ group. We do not know of two abstractly commensurable
groups, one of which is $\CATZ$, and the other of which is not.  We
promote the following question.
\begin{question} \label{TheQuestion2}
Is the property of being a $\CATZ$ group an invariant of the
abstract commensurability class of a group?
\end{question}
Some relevant results in the literature show that two natural
approaches to this question do not work in general.  If G acts
geometrically on a $\CATZ$ space $X$ and $G'$ is a finite extension
with $[G':G] = n$, then $G'$ acts properly and isometrically on the
$\CATZ$ space $X^n$ with the product metric \cite{SerreConstruction}
\cite[p.190]{KenBrown}. However, proving this action is cocompact is
either difficult or impossible in general. In
\cite{BestvinaKleinerSageev}, the authors give examples of the
following type: $G$ is a group acting faithfully and geometrically
on a $\CATZ$ space $X$, $G'$ is a finite extension of $G$, yet $G'$
does not act faithfully and geometrically on $X$. However, $G'$ may
act faithfully and geometrically on some other $\CATZ$ space.

\end{rem}
\begin{rem}\label{IsomorphismRemark}
The fact that $\Aut(F_2)$ is isomorphic to $\Aut(W_3)$ appears to be
well-known in certain mathematical circles, but is rarely recorded
explicitly. We now outline a proof: the subgroup $E \leq W_3$ of
even length elements is isomorphic to $F_2$, characteristic in $W_3$
and $C_{W_3}(E) = \{1\}$; it follows from \cite[Lemma 1.1]{Rose}
that the induced homomorphism $\pi\co\Aut(W_3) \to \Aut(E)$ is
injective; one easily confirms that the image of $\pi$ contains a
set of generators for $\Aut(E)$, and hence $\pi$ is an isomorphism.
A topological proof may also be constructed using the fact that the
subgroup $E$ of even length words in $W_3$ corresponds to the 2-fold
orbifold cover of the the orbifold $S^2(2, 2, 2, \infty)$ by the
once-punctured torus.
\end{rem}

The authors would like to thank Jason Behrstock and Martin Bridson
for pointing out the examples in \cite{KapovichLeeb} and
\cite[p.258]{MartinsBook}\cite{BridsonAlonso} and Luisa Paoluzzi for discussions regarding \cite{CrispPaoluzzi}.

\section{$\Aut(B_4)$ is a $\CATZ$ group}\label{AutB4Section}


We shall describe an apt presentation of $B_4$, give a concise
combinatorial description of Brady's space $X_0$, describe the
faithful geometric action of $\Inn(B_4)$ on $X_0$ and, finally,
introduce an isometry of $X_0$ to extend the action of $\Inn(B_4)$
to a faithful geometric action of $\Aut(B_4)$.

The interested reader will find an informative, and rather more
geometric, account of $X_0$ and the associated action of $\Inn(B_4)$
in \cite{CrispPaoluzzi}.

\subheading{An apt presentation of $B_4$} A standard presentation of
the group $B_4$ is
\begin{equation}\label{PresentationOfB4}
\langle a, b, c \; | \; aba = bab, bcb = cbc, ac = ca\rangle.
\end{equation}
Introducing generators $d = (ac)^{-1} b (ac)$, $e = a^{-1}ba$ and $f
= c^{-1}bc$, one may verify that $B_4$ is also presented by
\begin{equation}\label{PresentationOfG}
\begin{split}
\langle a, b, c, d, e, f  \; | \; & ba = ae=eb, de=ec =cd, bc = cf = fb, \\
 & df = fa = ad, ca = ac, ef = fe \rangle.
\end{split}
\end{equation}
We set $x = bac$ and write $\X \subset B_4$ for the
infinite cyclic subgroup generated by $x$.  The center of $B_4$ is
the infinite cyclic subgroup generated by $x^4$.

\subheading{The space $X_0$} Consider the 2-dimensional piecewise
Euclidean CW-complex $X_0$ constructed as follows:
\begin{itemize}
\item [(0-S)] the vertices of $X_0$ are in one-to-one correspondence with the left cosets of $\X$ in $B_4$---
we write $v_{g\X}$ for the vertex corresponding to the coset $g\X$;
\item [(1-S)] distinct vertices $v_{g_1\X}$ and $v_{g_2\X}$ are connected by an edge of unit length if and
only if there exists an element $\ell \in \{a, b, c, d, e, f\}^{\pm
1}$ such that $g_2^{-1} g_1 \ell \in \X$;
\item [(2-S)] three vertices $v_{g_1\X}$, $v_{g_2\X}$ and $v_{g_3\X}$ are the vertices of a
Euclidean (equilateral) triangle if and only if the vertices are
pairwise adjacent.
\end{itemize}

\begin{figure} \centering
\includegraphics[scale=0.65]{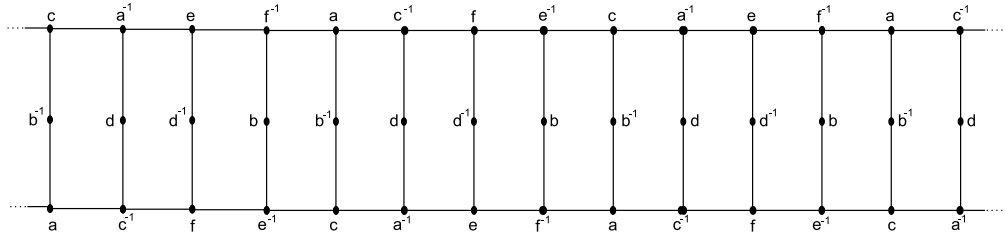}
\caption{A covering of the link of $v_\X$ in $X_0$.
\label{LinkOpenedUpDiagram}}
\end{figure}

The link of the vertex $v_{\X}$ in $X_0$, just like the link of each
vertex in $X_0$, consists of twelve vertices (one for each of the
cosets represented by elements in $\{a, b, c, d, e, f\}^{\pm 1}$)
and sixteen edges (one for each of the distinct ways to spell $x$ as
a word of length three in the alphabet $\{a, b, c, d, e, f\}$---see
\cite{CrispPaoluzzi} for more details). It can be viewed as the
1-skeleton of a M\"obius strip. In Figure \ref{LinkOpenedUpDiagram}
we depict the infinite cyclic cover of the link of $v_\X$. Each
vertex with label $g$ in the figure lies above the vertex $v_{g\X}$
in the link of $v_\X$. The link is formed by identifying identically
labeled vertices and identifying edges with the same start and end
points.

That $X_0$ is $\CATZ$ follows most naturally from the alternative
construction of $X_0$ described in detail in \cite{CrispPaoluzzi}.
Alternatively, a complex constructed from isometric Euclidean triangles is $\CATZ$ if and only if it is simply-connected and satisfies the `link condition' \cite[Theorem II.5.4, pp.206]{MartinsBook}. For a 2-dimensional complex, the link condition requires that each injective loop in the link of a vertex has length at least $2\pi$, where edges in a link are assigned the length of the angle they subtend \cite[Lemma II.5.6, pp.207]{MartinsBook}. It is easily seen that $X_0$ satisfies the link condition because each injective loop in Figure \ref{LinkOpenedUpDiagram} crosses at least 6 edges and each edge has length $\pi/3$.  Thus one might show that $X_0$ is $\CATZ$
by showing that it is simply-connected.  We shall not digress from
the task at hand to provide such an argument.

\subheading{Brady's faithful geometric action of $\Inn(B_4)$ on
$X_0$}  We shall describe Brady's faithful geometric action of
$\Inn(B_4)$ on $X_0$.  We shall do so by describing an isometric
action $\rho\co B_4 \to \Isom(X_0)$ such that the image of $\rho$ is
a properly discontinuous and cocompact subgroup of $\Isom(X_0)$
which is isomorphic to $\Inn(B_4)$.

It follows immediately from (1-S) that, for each $g \in B_4$, the
``left-multipli-cation by $g$'' map on the 0-skeleton of $X_0$,
$g_1\X \mapsto g g_1\X$, extends to a simplicial isometry of the
1-skeleton of $X_0$. It follows immediately from (2-S) that any
simplicial isometry of the 1-skeleton of $X_0$ extends to a
simplicial isometry of $X_0$. We write $\phi_g$ for the isometry of
$X_0$ determined by $g$ in this way, and we write $\rho\co B_4 \to
\Isom(X_0)$ for the map $g \mapsto \phi_g$.  We compute that
$\rho(g_1 g_2)(v_{g\X}) = v_{g_1 g_2 g\X} = \rho(g_1)
\rho(g_2)(v_{g\X})$ for each $g_1, g_2, g  \in B_4$, so $\rho$ is a
homomorphism. Further, $\phi_g(v_\X) = v_{g\X}$ for each $g  \in
B_4$, so the vertices of $X_0$ are contained in a single
$\rho$-orbit. It follows that $\rho$ is a cocompact isometric action
of $B_4$ on $X_0$.

\begin{figure} \centering
\includegraphics[scale=0.65]{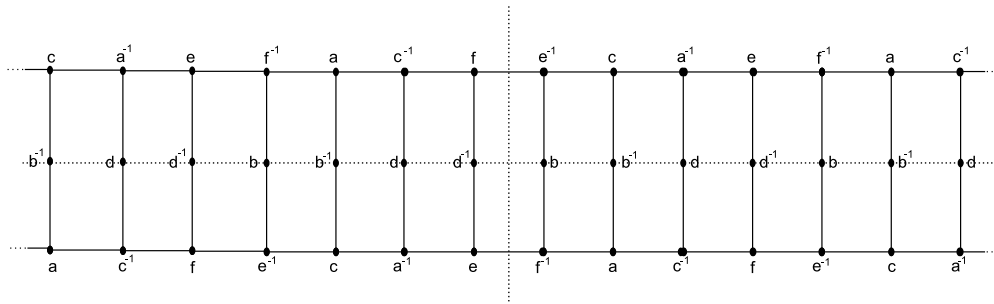}
\caption{A covering of the link of the vertex $v_\X$ and the fixed
point sets of some reflections. \label{LinkOpenedUpWithLineDiagram}}
\end{figure}

To show that the image of $\rho$ is isomorphic to $\Inn(B_4)$, it
suffices to show that the kernel of $\rho$ is exactly the center of
$B_4$. One easily computes that $\rho(x^4)$ is the identity isometry
of $X_0$. Thus the kernel of $\rho$ contains the center of $B_4$. It
is also clear that the stabilizer of $v_{\X}$, which contains the
kernel of $\rho$, is the infinite subgroup $\X $. So to establish
that the kernel of $\rho$ is exactly the center of $B_4$, it
suffices to show that $\phi_x, \phi_{x^2}$ and $\phi_{x^3}$ are
non-trivial and distinct isometries of $X_0$. We achieve this by
showing that these elements act non-trivially and distinctly on the
link of $v_{\X}$ in $X_0$.  We compute that $x$ acts as follows on
the cosets corresponding to vertices in the link of $v_\X$, where
$\delta = \pm 1$:
$$a^\delta \X \mapsto e^\delta \X \mapsto c^\delta \X \mapsto f^\delta \X \mapsto a^\delta \X \hbox{ and } b^\delta \X \leftrightarrow d^\delta \X.$$
Thus the restriction of $\phi_x$ to the link of $v_\X$ may be
understood, with reference to Figure
\ref{LinkOpenedUpWithLineDiagram}, as translation two units to the
right followed by reflection across the horizontal dotted line.
It follows that $\phi_x, \phi_{x^2}, \phi_{x^3}$ are non-trivial and
distinct isometries of $X_0$, as required

We next show that the image of $\rho$ is a properly discontinuous
subgroup of $\Isom(X_0)$. Now, the action $\rho$ is not properly
discontinuous because, as noted above, the $\rho$-stabilizer of
$v_{\X}$ is the infinite subgroup $\X$ (so infinitely many elements
of $B_4$ fail to move the unit ball about $v_{\X}$ off itself).  But
the image of $\X$ under the map $B_4 \to \Inn(B_4)$ has order four.
It follows that the image of $\rho$ is a properly discontinuous
subgroup of $\Isom(X)$.

Thus we have that the image of $\rho$ is a properly discontinuous
and cocompact subgroup of $\Isom(X_0)$ which is isomorphic to
$\Inn(B_4)$.

\subheading{Extending $\rho$ by finding one more isometry} It was
shown in \cite{AutomorphismsOfBraidGroups} that the unique
non-trivial outer automorphism of $B_n$ is represented by the
automorphism which inverts each of the generators in Presentation
 (\ref{PresentationOfB4}). Consider the automorphism $\tau \in \Aut(B_4)$ determined by
$$a \mapsto a^{-1}, \;\;\; b \mapsto d^{-1}, \;\;\; c \mapsto
c^{-1}, \;\;\; d \mapsto b^{-1}, \;\;\; e \mapsto f^{-1}, \;\;\; f
\mapsto e^{-1}.$$  Note that $\tau$ is achieved by first applying
the automorphism which inverts each of the generators $a$, $b$ and
$c$ and then applying the inner automorphism $w \mapsto (ac)^{-1} w
(ac)$ for each $w \in B_4$.  It follows that $\tau$ is an involution
which represents the unique non-trivial outer automorphism of $B_4$.
Writing $J := B_4 \rtimes_\tau \Integer_2$, we have $\Aut(B_4) \cong
J/\langle x^4 \rangle$.  We identify $B_4$ with its image in $J$.

The automorphism $\tau \in \Aut(B_4)$ permutes the elements of $\{a,
b, c, d, e, f\}^{\pm 1}$ and maps the subgroup $\X$ to itself (in
fact, $\tau(x) = x^{-1}$).  It follows from (1-S) that the map
$v_{g_1\X} \mapsto v_{\tau(g_1)\X}$ on the 0-skeleton of $X_0$
extends to a simplicial isometry of the 1-skeleton of $X_0$, and
hence also to a simplicial isometry $\theta$ of $X_0$.  We compute
that $\theta \phi_{g} \theta = \phi_{\tau(g)}$ for each $g \in B_4$.
Thus we have an isometric action $\rho'\co J \to \Isom(X_0)$ given
by $$g \mapsto \phi_g \hbox{ for each } g \in B_4, \hbox{ and } \tau
\mapsto \theta.$$ We also compute that the restriction of $\theta$
to the link of $v_\X$ may be understood as reflection across the
vertical dotted line shown in Figure
\ref{LinkOpenedUpWithLineDiagram}. It follows that $\theta$ is a
non-trivial isometry of $X_0$ which is distinct from $\phi_x,
\phi_{x^2}$ and $\phi_{x^3}$. Thus the kernel of $\rho'$ is still
the center of $B_4$, and the image of $\rho'$ is a properly
discontinuous and cocompact subgroup of $\Isom(X_0)$ which is
isomorphic to $\Aut(B_4)$. Hence we have a faithful geometric action
of $\Aut(B_4)$ on $X_0$, as required.

\bibliographystyle{abbrv}
\bibliography{AutF2IsCAT0Bib}
\end{document}